\documentclass{article}
\usepackage{amsmath, epsfig,amssymb, multicol}
\newtheorem{lemma}{Lemma}
\newtheorem{theorem}{Theorem}
\newtheorem{conjecture}{Conjecture}

\newtheorem{corollary}{Corollary}
\newtheorem{property}{Property}

\newtheorem{claim}{Claim}

\newcommand{\Z}{{\mathbb Z}}

\newcommand{\E}{{\mathbb E}}

\title{Monochromatic 4-term arithmetic progressions \linebreak in 2-colorings of $\mathbb Z_n$}
\author{Linyuan Lu
\thanks{University of South Carolina, Columbia, SC 29208,
({\tt lu@math.sc.edu}). This author was supported in part by NSF
grant DMS 1000475. }
  \and Xing Peng
\thanks{University of South Carolina, Columbia, SC 29208,
({\tt pengx@mailbox.sc.edu}).This author was supported in part by
NSF grant  DMS 1000475. } }
\begin{document}
\maketitle
\begin{abstract}
  This paper is motivated by a recent result of Wolf \cite{wolf} on
  the minimum number of monochromatic 4-term arithmetic progressions
  (4-APs, for short) in $\Z_p$, where $p$ is a prime number.  Wolf
  proved that there is a 2-coloring of $\Z_p$ with $0.000386\%$
  fewer monochromatic 4-APs than random $2$-colorings; 
 the proof is probabilistic and non-constructive.
  In this paper, we present an explicit and simple construction
  of a $2$-coloring with $9.3\%$ fewer monochromatic 4-APs
  than random $2$-colorings.  This problem leads us to consider the minimum
  number of monochromatic 4-APs in $\Z_n$ for general $n$.  We obtain both lower
  bound and upper bound on the minimum number of monochromatic 4-APs in all 2-colorings of $\Z_n$.
 Wolf proved  that any $2$-coloring of $\Z_p$
  has at least  $(1/16+o(1))p^2$  monochromatic 4-APs. We improve this lower bound
 into  $(7/96+o(1))p^2$.

 Our results on $\Z_n$ naturally apply to the similar problem on $[n]$ (i.e., $\{1,2,\ldots, n\}$).  In 2008,
 Parillo, Robertson, and Saracino \cite{prs} constructed a
 $2$-coloring of $[n]$ with 14.6\% fewer monochromatic 3-APs
 than random $2$-colorings. In 2010, Butler, Costello, and Graham
 \cite{BCG} extended their methods and used an extensive computer search
 to construct a $2$-coloring of $[n]$ with 17.35\% fewer monochromatic
 4-APs (and 26.8\% fewer monochromatic 5-APs) than random
 $2$-colorings. Our construction gives a $2$-coloring of $[n]$ with 33.33\%
 fewer monochromatic 4-APs (and 57.89\% fewer monochromatic 5-APs) than random 2-colorings.

\end{abstract}

\section{Introduction}

Let $G$ be a finite subset of a commutative group.  For any integer
$k\geq 3$, a {\em $k$-term arithmetic progression} (or $k$-AP, for
short) is an (ordered) sequence of  $k$ elements in $G$ of the form $(a, a+d, \ldots,
a+(k-1)d)$, where $a$ is the {\em first element} and
$d$ is the {\em common difference}.
 A $2$-coloring of $G$ is a map $c\colon G\to \{0,1\}$.
A $k$-AP $(a, a+d, \ldots, a+(k-1)d)$ is {\em monochromatic} if
$c(a)=c(a+d)=\cdots=c(a+(k-1)d)$.  Let $m_k(G,c)$ be the number of
monochromatic $k$-APs in the 2-coloring $c$.  A natural question is
how small $m_k(G,c)$ can be?  Let $AP_k(G)$ be the number of all
$k$-APs in $G$. Define
\begin{equation}
  \label{eq:mk}
m_k(G):=\min_{c}\frac{  m_k(G,c)}{AP_k(G)}.
\end{equation}
We are interested in the asymptotic value of $m_k(G)$ as $|G|$
approaches infinity.  (This is similar to those questions on  Schur
Triples \cite{da, RZ, thanatipanonda} or on general patterns
\cite{BCG, fgr}.)

In this paper, we consider only the cases that $G=[n]$ and $G=\Z_n$.
Here $[n]=\{1,2,\ldots,n\}$ and $\Z_n$ is the cyclic group of order
$n$. When $n$ is a prime number $p$, we  write $\Z_n$ as $\Z_p$ for
emphasis.  A $k$-AP is {\em degenerated} if it contains repeated
terms; it is {\em non-degenerated} otherwise.  The {\em mirror
  image} of a $k$-AP $(a,a+d,\ldots, a+(k-1)d)$ is another $k$-AP
$(a+(k-1)d,\ldots, a+d,a)$.  Here we allow $k$-APs to be degenerated;
 a $k$-AP differs from its mirror image
except for $d=0$.  In contrast,  many papers require $k$-APs to be
non-degenerated and treat a $k$-AP  the same as its mirror image.
The two different definitions of $k$-APs derive two different versions of
$m_k(G)$.  However, they are asymptotically equivalent as $|G|$ goes
to infinity;  this is because the number of degenerated $k$-APs is
only $O(n)$ while the number of all APs is $\Omega(n^2)$.  A $k$-AP
$(a, a+d, \ldots, a+(k-1)d)$ is parametrized by a pair $(a,d)$.  The
parameter space of all $k$-APs in $[n]$ can be described as
$\{(a,d)\colon 1\leq a \leq n,  1\leq a+(k-1)d\leq n\}.$ A $k$-AP
$(a, a+d, \ldots, a+(k-1)d)$ in $[n]$ is degenerated if and only if
$d=0$. The parameter space of all $k$-APs in $\Z_n$ is simply 
$\Z_n^2$. A $k$-AP $(a, a+d, \ldots, a+(k-1)d)$ in $\Z_n$ is
degenerated if $jd\equiv 0 \mod n$ for some $0 \leq j \leq k-1$. In
both cases, the number of degenerated $k$-APs is $O(n)$.

 Random $2$-colorings of $[n]$ (or $\Z_n$) give the following  upper bounds.
\begin{eqnarray}
  \label{eq:randn}
  m_k([n])&\leq& 2^{1-k}+o(1);\\
  \label{eq:randzn}
  m_k(\Z_n)&\leq& 2^{1-k}+o(1).
\end{eqnarray}

Van der Waerden's number \cite{waerden} $W=W(2,k)$ can be used to
provide a lower bound on $m_k([n])$. For example, using a double
counting method, one can prove $m_k([n])\geq
\frac{2(k-1)}{W^3}+o(1)$ (see \cite{BCG}). A similar argument can
show $m_k(\Z_n)\geq \frac{2(k-1)}{W^2}+o(1)$. These bounds are
usually too weak; stronger bounds exist for $k=3$ and $k=4$.



The case $\Z_p$  is of particular interest.  The number of
monochromatic $3$-APs in $\Z_p$ depends only on the size of the
coloring classes, but not on the coloring itself (see \cite{da}).
Namely, if $c$ is a 2-coloring of $\Z_p$ such that the size of red
class is $\alpha p$, then we have
\begin{equation}
  \label{eq:m3p}
m_3(\Z_p,c)=(1-3\alpha+3\alpha^2) p^2.
\end{equation}
The minimum is achieved when $\alpha$ is closed to $\frac{1}{2}$.
Thus $m_3(\Z_p)$ is achieved by random 2-colorings.

For $k=4$, Wolf \cite{wolf} proved that for any sufficiently large
prime number $p$, we have
\begin{equation}
  \label{eq:wolf}
\frac{1}{16}+o(1) \leq m_4(\Z_p)\leq
\frac{1}{8}(1-\frac{1}{259200})+o(1).
\end{equation}
This lower bound improved a previous lower bound
due to Cameron, Cilleruelo, and Serra \cite{cameron},
\begin{equation}
  \label{eq:cameron} m_4(\Z_n) \geq \frac{2}{33}+o(1),
\end{equation}
 where $n$ is relatively prime to $6$ and large enough. (Cameron, Cilleruelo, and
Serra's result actually holds for any Abelian group of order $n$
provided $\mbox{gcd}(n,6)=1$.)

Wolf's upper bound indicates that $m_4(\Z_p)$ is not achieved by
random 2-colorings.  This is a nice result. 
However, the quantity is only slightly less than
$\frac{1}{8}$
--- the   density of monochromatic $4$-APs in random $2$-colorings. Her
method for the upper bound relies heavily on  the method initialized
by Gowers  (see \cite{wolf}). The existence of
such $2$-coloring is proved  by probabilistic methods; it is
non-constructive.

To get a better upper bound for $m_k(\Z_n)$, we introduce a
construction consisting of periodic blocks. For a fixed $b$, let $B$
be a good $2$-coloring of $\Z_b$ with $m_k(\Z_b)b^2$ monochromatic
$k$-APs. (Here $B$ is viewed as a 0-1 string of length $b$.)

Write $n=bt+r$ with $0\leq r\leq
b-1$. We consider the following periodic construction $c$
\begin{equation}
  \label{eq:pcon}
\underbrace{BB\cdots B}_tR.
\end{equation}
Here $R$ is any bit-string of length $r$.

If $n$ is divisible by $b$, then $R$ is empty. In this
case, it is easy to see that the periodic construction above gives
$m_k(\Z_b)n^2$ monochromatic $k$-APs. Thus, we have
\begin{equation}
  \label{eq:divisible}
m_k(\Z_n)\leq m_k(\Z_b) \quad \mbox{ if }b\mid n.
\end{equation}

If $n$ is not divisible by $b$, then the computation of
$m_k(\Z_n,c)$ is more complicated in general. Note that the number
of $k$-APs containing some element(s) in $R$ is a lower order term
as $n$ goes to infinity; the value $m_k(\Z_n,c)$ can be still
determined asymptotically by $B$. (See Lemma \ref{upper} and
\ref{le3}.)

Two colorings $c$ and $c'$ of $\Z_n$ are {\em isomorphic} if there
is an integer $m$  such that $\mbox{gcd}(m,n)=1$ and  $c'(v)=c(mv)$
for any $v\in \Z_n$. Two colorings $c$ and $c'$ of $\Z_n$ are {\em
conjugated} if $c'(v)=1-c(v)$ for any $v\in \Z_n$. It is clear that
$m_k(\Z_n,c)=m_k(\Z_n,c')$, whenever $c$ and $c'$ are isomorphic or
conjugated to each other. To find a good coloring $B$, we implement
an efficient bread-first search algorithm reducing isomorphic
copies. From the proof of the lower bound for $m_4(\Z_n)$, we need
pay attention to those $n$ divisible by $4$. Using this efficient
program, we find
 a  good $2$-coloring $B_{20}$ of $\Z_{20}$ for 4-APs,
$$B_{20}=(1,1,1,0,1,1,0,1,1,1,0,0,0,1,0,0,1,0,0,0).$$
This coloring $B_{20}$ gives $m_4(\Z_{20})=\frac{9}{100}.$

We also search the coloring of $\Z_p$ without any non-degenerated
monochromatic $4$-APs. At $p=11$, there is a unique coloring with
this property up to isomorphisms. Since $0$'s and $1$'s are not
balanced in this coloring, we search  good colorings in  $\Z_{22}$ instead. We
found a good $2$-coloring $B_{22}$ of $\Z_{22}$ for 4-APs,
$$
B_{22}=(1,1,1,0,1,1,0,1,0,0,0,1,1,1,0,1,0,0,1,0,0,0).
$$
In this coloring, all monochromatic $4$-APs in $B_{22}$ are degenerated;
there are  22 monochromatic $4$-APs with $d=0$ and 20 monochromatic
 $4$-APs with $d=11$. This
coloring $B_{22}$ gives $m_4(\Z_{22})=\frac{42}{22^2}=\frac{21}{242}.$

The following theorem improves both Wolf's lower bound and upper
bound on $m_4(\Z_p)$. Our lower bound is obtained by combining
Wolf's elegant method and an  exhaustive search. Our upper bound is
proved by a novel method of analyzing the number of monochromatic
$k$-APs in the periodic construction (\ref{eq:pcon}).

\begin{theorem}\label{prime}
If $p$ is  prime and large enough, then we have
  \begin{equation}
    \label{eq:m4p}
0.07291666<    \frac{7}{96} \leq
m_4(\Z_p)\leq\frac{17}{150}+o(1)<0.1133334.
  \end{equation}
\end{theorem}

In fact, our methods naturally lead (asymptotic) bounds on $m_4(\Z_n)$ for general $n$.
The results depend on
$n$ case-wisely. For simplicity, we split it into two theorems: one
on the lower bound and the other one on the upper bound.

\begin{theorem}\label{t:m4lb}
If $n$ is sufficiently large, then we have
$$m_4 (\Z_n)\geq \left\{
  \begin{array}[c]{ll}
   \frac{7}{96} & \mbox{ if $n$ is not divisible by 4},\\

\vspace{-0.3cm}&\\

    \frac{2}{33} & \mbox{ if $n$ is divisible by 4.} \\
  \end{array}
\right.
$$
\end{theorem}

Here is a theorem  for the upper bound on $m_4(\Z_n)$ for general
$n$.

\begin{theorem}\label{t:m4ub}
For $n$ sufficiently large, we have
$$m_4 (\Z_n)\leq \left\{
  \begin{array}[c]{ll}
    \frac{17}{150}+o(1)< 0.1133334& \mbox{ if  $n$ is odd,}  \\
    \vspace{-0.3cm}&\\
  \frac{8543}{72600}+o(1)<0.1176722 & \mbox{ if $n$ is even.}
  \end{array}
\right.
$$
\end{theorem}
Note that Theorem \ref{prime} is a corollary of Theorem \ref{t:m4lb}
and \ref{t:m4ub}. The upper bound above is small enough to beat the
bound $\frac{1}{8}$
 reached by random 2-colorings. Using inequality (\ref{eq:divisible}), we can
get a much better bound for certain $n$'s. For example,
$$m_4(\Z_n) \leq
\left\{
  \begin{array}[c]{ll}
    0.09 & \mbox{ if  } 20\mid n,  \\
    0.086777 & \mbox{ if  } 22\mid n,  \\
  \end{array}
\right.
$$

For $m_5(\Z_n)$, we use the periodic construction with
the following good coloring  of $\Z_{74}$:\\
{\small $B_{74}=(1,1,1,1,0,1,1,1,0,0,0,0,1,0,1,1,0,0,1,0,1,0,1,0,0,1,1,0,1,0,0,0,0,1,1,1,0,$\\
 $\hspace*{1.1cm}  1,1,1,1,0,1,1,1,0,0,0,0,1,0,1,1,0,0,1,0,0,0,1,0,0,1,1,0,1,0,0,0,0,1,1,1,0).$}

All monochromatic $5$-APs  in $B_{74}$ are
degenerated ones. Among them there are  74 $5$-APs with $d=0$ and 72 $5$-APs with
$d=37$. This coloring gives
$m_5(\Z_{74})=\frac{146}{74^2}=\frac{73}{2738}$. We have  the
following theorem.





\begin{theorem} \label{t:m5ub}
If $n$ is  sufficiently large, then we have
$$m_5 (\Z_n)\leq \left\{
  \begin{array}[c]{ll}
  \frac{3629}{65712}+o(1)<0.055226 & \mbox{ if  $n$ is odd,}  \\
   \vspace{-0.3cm}&\\
  \frac{3647}{65712}++o(1)<0.0554998& \mbox{ if $n$ is even.}
  \end{array}
\right.
$$
\end{theorem}


Once again, the upper bound above is small enough to beat the bound $\frac{1}{16}$,
which is  reached by random 2-colorings.
Using inequality (\ref{eq:divisible}),
we can get a much better upper bound for certain $n$. For example, we have
$$m_5(\Z_n)\leq \frac{73}{2738}=0.026661\cdots \mbox{ if } 74 \mid n.$$
The following theorem  gives the best lower bounds (for some $n$'s).
 \begin{theorem} \label{tzb}
We have
\begin{eqnarray*}
   \varliminf_{n\to \infty}m_4(\Z_n)\leq \frac{1}{12}\\
   \varliminf_{n\to \infty}m_5(\Z_n)\leq \frac{1}{38}.
\end{eqnarray*}
 \end{theorem}
In fact, we show for any $\epsilon$, there is an odd integer $n$ with
$m_4(\Z_n)\leq \frac{1}{12}+\epsilon$.  Combining this result with
Theorem \ref{t:m4lb}, we get
$$\frac{7}{96}  \leq \inf\{m_4(\Z_n)\colon n \mbox{ is not divisible by 4 }\}\leq \frac{1}{12}.$$
Note the gap is pretty small. Here we conjecture that the upper bound is
tight.
\begin{conjecture}
$\inf\{m_4(\Z_n)\colon n \mbox{ is not divisible by 4 }\}= \frac{1}{12}.$
\end{conjecture}
Maybe it is true even if the condition that ``$n$ is not divisible by 4'' is removed.

The periodic construction also works for $m_k([n])$. We have
\begin{lemma} \label{l0}
For any $k\geq 3$ and any positive integer $b$, we have
$$\varlimsup_{n\to\infty}m_k([n])\leq m_k(\Z_b).$$
In particular, we have
$$\varlimsup_{n\to\infty}m_k([n])\leq \varliminf_{n\to\infty} m_k(\Z_n).$$
\end{lemma}

When we consider the similar problems for $[n]$, the $k$-AP and its
mirror image are often not distinguished in the literature. To avoid
the ambiguity, we call a $k$-AP (in $[n]$) with $d>0$
 an {\em increasing} $k$-AP. For $k\geq 3$, let $c_k$ be the largest
number satisfying ``for any $\epsilon>0$, there is a sufficiently large $n$
such that  any $2$-coloring of $[n]$ contains at least
$(c_k-\epsilon)n^2$ monochromatic increasing $k$-APs''. 
Since $[n]$ has $(\frac{1}{2(k-1)}+o(1))n^2$ increasing
$k$-APs, it is equivalent to say
\begin{equation}
  \label{eq:ck}
  c_k=\frac{1}{2(k-1)}\varlimsup_{n\to\infty} m_k([n]).
\end{equation}

In 2008, Parillo, Robertson, and  Saracino \cite{prs}
proved
\begin{equation}
  \label{eq:c3}
0.05111<\frac{1675}{32768} \leq c_3\leq \frac{117}{2192}<0.053376.
\end{equation}
Their construction was generalized by Butler, Costello, and
Graham, who \cite{BCG} proved $c_4<  0.0172202\ldots$ and
$c_5<0.005719619\ldots$ via an extensive computation. Both bounds
beat
 random $2$-colorings.

Combining Theorem \ref{tzb} with Lemma \ref{l0},  we have
 \begin{eqnarray}
   \label{eq:c4}
   c_4&\leq& \frac{1}{72}=0.01388888\ldots,\\
   c_5 &\leq& \frac{1}{304}=0.003289474\ldots.
 \end{eqnarray}
These numerical results indicate that the periodic
construction is often better than the block construction used in \cite{BCG}. 
We believe the following conjecture holds.
\begin{conjecture}
For fixed $k\geq4$,  we have  $\varlimsup_{n\to\infty}m_k([n])=
\varliminf_{n\to\infty} m_k(\Z_n)$.
\end{conjecture}


Bounding $c_3$ is very different from bounding $c_4$.
This conjecture above is not true for $k=3$. We have
the following theorem.
\begin{theorem}\label{t:z3}
If the integer $n$ is large enough, then  any $2$-coloring of $\Z_n$
contains at least $\frac{1}{4}n^2$ monochromatic arithmetic
progressions. In particular, we have
\begin{equation}
  \label{eq:m3}
m_3(\Z_n)=\frac{1}{4}+o(1).
  \end{equation}
\end{theorem}

With the help of computer search, we found three good 2-colorings
$B_{20}$, $B_{22}$, and $B_{74}$,  which are used  as building
blocks in constructing good 2-colorings of $\Z_n$ and $[n]$. The data
in Table \ref{tab:3}, \ref{tab:4}, \ref{tab:5}, and \ref{tab:6},  can be easily
verified by anyone with limited programming experience. Some lower
bound requires nontrivial exhaustive search in the same way as
Cameron, Cilleruelo, and Serra \cite{cameron} proved the previous
lower bounds. However, those lower bounds using an exhaustive
computer search are not the focus of this paper.

The organization of the paper is following. In section 2, we will
prove a necessary lemma and Theorem \ref{t:z3}.
 In the section 3, we
first prove a lemma and a corollary on counting lattice points in a
polygon; then we prove Theorem \ref{t:m4ub} for odd $n$ and  Theorem
\ref{t:m5ub}. In  section 4, we introduce a recursive construction
and then use it to prove Theorem \ref{tzb} and Theorem \ref{t:m4ub}
for even $n$. In the last section, we will deal with the lower
bounds and prove Theorem \ref{t:m4lb}.

\section{Notations and the proof of Theorem \ref{t:z3}}
Let $c\colon \mathbb Z_n \to \{0,1\}$ be a $2$-coloring of $\mathbb
Z_n$. The coloring $c$ is often viewed as a bit-string of length
$n$. For convenience, we say an element $v\in \Z_n$ is {\em red} if
$c(v)=0$ and  {\em blue} if $c(v)=1$. The coloring $c$ induces a
partition  $\Z_n=A\cup B$, where $A$ is the set of red elements
while $B$ is the set of blue elements.


Let $k \geq 3$ be an integer  and $|A|=\alpha n$. We have
$|B|=(1-\alpha)n$.

For each $1\leq i \leq k$,  let $A_i$ (or $B_i$) be the set of all
$k$-APs  whose   $i$-th number  is red (or blue), respectively; we
have
\begin{eqnarray}
\label{eq:Ai}
|A_i|&=&\alpha n^2, \\
\label{eq:Bi} |B_i|&=&(1-\alpha) n^2.
\end{eqnarray}
\begin{lemma}\label{l1}
  For $1\leq i <j \leq k$, if $\mbox{gcd}(j-i,n)=1$, then
we have
\begin{eqnarray}
\label{eq:Aij}
  |A_i\cap A_j|&=& \alpha^2 n^2,\\
\label{eq:Bij}
  |B_i\cap B_j|&=& (1-\alpha)^2 n^2.
\end{eqnarray}
If $\mbox{gcd}(j-i,n)\not=1$, then we have
\begin{eqnarray}
\label{eq:Aij1}
  |A_i\cap A_j|&\geq& \alpha^2 n^2,\\
\label{eq:Bij1}
  |B_i\cap B_j|&\geq&(1-\alpha)^2 n^2.
\end{eqnarray}
\end{lemma}
{\bf Proof:} For $1\leq i <j \leq k$,  the value of  $|A_i\cap A_j|$
equals the number of $k$-APs whose $i$-th and $j$-th terms are red.
If $\mbox{gcd}(j-i,n)=1$, then every ordered pair of elements
(distinct or not)
 in $\Z_n$ can be extended into a unique $k$-AP whose $i$-th and
$j$-th terms are the given pair. Note the number of  ordered pairs
of red (and blue)  elements  is exactly $\alpha^2 n^2$ (and
$(1-\alpha)^2 n^2$), respectively. Equations (\ref{eq:Aij}) and
(\ref{eq:Bij}) follow.

 If $\mbox{gcd}(j-i,n) \not=1$, then  every pair of elements in $\Z_n$
may or may not be extended into a $k$-AP whose $i$-th and $j$-th
terms are the given pair.  Let $r=\mbox{gcd}(j-i, n)$. For $0 \leq
l\leq r-1$, let $x_l$ be the number elements $z$ in $\Z_n$ such that
$z$ is red and $z\equiv l \mod r$.  For any pair $(z_1,z_2)$, the
elements $z_1$ and $z_2$ are the $i$-th and $j$-th elements of an
arithmetic progression  if
\begin{equation}
  \label{eq:d}
z_2-z_1=(j-i)d,
\end{equation}
for some element $d$ in $\Z_n$. Equivalently, $z_2-z_1 \equiv 0 \mod
r$. Moreover, if $z_2-z_1\equiv 0 \mod r$, then equation
(\ref{eq:d}) has $r$ solutions. It follows that
\begin{eqnarray}
  \label{eq:Aij2}
  |A_i\cap A_j|&=& r\sum_{l=0}^{r-1}x_l^2\\
\nonumber
&\geq& (\sum_{l=0}^{r-1}x_l)^2\\
\nonumber
&=& \alpha^2 n^2.
\end{eqnarray}
Equation (\ref{eq:Bij1}) can be proved similarly.
\hfill $\square$

\noindent
{\bf Proof of theorem \ref{t:z3}:}  Observe that if we  assign red
and blue to each number equally likely,  then the expected value of
$m_3(\Z_n,c)$ is $n^2/4 +O(n)$. Therefore, there is a 2-coloring $c$
such that $m_3(\Z_n,c) \leq n^2/4+O(n)$, that is  $m_3(\Z_n) \leq 1/4+O(1/n)$.

For the other direction,  let $c$ be any 2-coloring of $\Z_n$. We
use the notations $\alpha$, $A_i$, and $B_i$ defined in the
beginning of this section.

 We have the following inclusion-exclusion formula.
\begin{equation} \label{eq:1}
|A_1 \cup A_2 \cup A_3|= \sum_{i=1}^3|A_i|-\sum_{1 \leq i < j \leq
3} |A_i \cap A_j|+|A_1 \cap A_2 \cap A_3|.
\end{equation}
Note that $A_1 \cup A_2 \cup A_3 = \overline{B_1\cap B_2\cap B_3}$
and $\left|\overline{B_1\cap B_2\cap B_3}\right|=n^2- |B_1\cap
B_2\cap B_3|.$ By the definition of $m_3(\Z_n,c)$, we have  $|A_1
\cap A_2 \cap A_3| + |B_1\cap B_2\cap B_3| = m_3(\Z_n,c)$. Applying
Lemma \ref{l1}, we have
\begin{eqnarray*}
m_3(\Z_n,c) &=& n^2 - \sum_{i=1}^3 |A_i|+\sum_{1 \leq i < j \leq
3} |A_i \cap A_j|\\
&\geq & n^2 - 3\alpha n^2 + 3 \alpha^2 n^2\\
&=& (1-3\alpha(1-\alpha))n^2.
\end{eqnarray*}
Note that $\alpha(1-\alpha)$ reaches the maximum value at
$\alpha=1/2$. We have $m_3(\Z_n,c)\geq n^2/4$. Therefore $m_3(\Z_n)
\geq 1/4$ and the lemma follows. \hfill $\square$


\section{Proofs of Theorem \ref{t:m4ub} and  Theorem \ref{t:m5ub}}

In this section, we will examine the number of  monochromatic $k$-APs
in the periodic construction (\ref{eq:pcon}).

\subsection{Proof of Lemma \ref{l0}}
We need a tool to count the grid points inside a polygon on the plane.

A point in ${\mathbb R}^2$ is a {\it  grid point} if both coordinates are integers.
Let $Q$ be a simple polygon whose vertices are grid points.
Let $A(Q)$ be the area of $Q$, $I(Q)$ be the number of grid points inside $Q$,
and $B(Q)$ be the number of grid points on the boundary of $Q$.
The classical Pick's theorem \cite{pick} states
$$A(Q)=I(Q)+\frac{B(Q)}{2}-1.$$

Intuitively, if $B(Q)$ is a lower order term, then $I(Q)\approx
A(Q)$. Let $P$ be a simple polygon in the plane ${\mathbb R}^2$. For
any $t>0$ and a point $v$, a new polygon $v+tP$ is  obtained by
first scaling $P$ by a factor of $t$ and then translating by a
vector $v$.  We have the following lemma.
\begin{lemma}
Suppose $P$ is a simple polygon with $m$ vertices and circumference
$L$.
 For any point $v$ and sufficiently large $t$, we have
$$|I(v+tP)-A(P)t^2|\leq 3Lt+ 5m.$$
\end{lemma}
{\bf Proof:} Since $P$ has $m$ vertices,  let $v_1,\ldots, v_m$ be
the vertices of the polygon $v+tP$. For $i=1,\ldots, m$, let $u_i$
be a grid point closest to $v_i$ (if there are more than one choice,
then break ties arbitrarily). We have $|u_iv_i|\leq
\frac{\sqrt{2}}{2}$. Let $Q$ be the polygon with vertices $u_1,
u_2,\ldots, u_m$. (For convenience, we write $v_{m+1}=v_0$ and
$u_{m+1}=u_0$.) The polygon $Q$ can be viewed as an approximation of
the polygon $v+tP$; thus $Q$ is simple for sufficiently large $t$.

Applying Pick's theorem to $Q$, we have
$$ A(Q)-I(Q)=\frac{B(Q)}{2}-1.$$
We observe that the number of  grid points on any  line segment of
length $l$ is at most $l+1$. We have
\begin{eqnarray*}
B(Q)&\leq& \sum_{i=1}^m(|u_iu_{i+1}|+ 1)  \\
&\leq& \sum_{i=1}^m(|v_iv_{i+1}|+ |u_iv_i| +|u_{i+1}v_{i+1}|+ 1)\\
&\leq& \sum_{i=1}^m(|v_iv_{i+1}| + \sqrt{2}+1)\\
&=& tL+(\sqrt{2}+1)m.
\end{eqnarray*}
Let $S_i$ be the convex region spanned by $v_i,v_{i+1}, u_i,
u_{i+1}$. Note $S_i$ is covered by four triangles $\Delta
u_iv_iv_{i+1}$, $\Delta u_iv_iu_{i+1}$,  $\Delta u_{i+1}v_{i+1}u_i
$, and $\Delta u_{i+1}v_{i+1}v_i$ exactly twice. We have
\begin{eqnarray*}
  A(S_i)&=&\frac{1}{2} (A(\Delta u_iv_iv_{i+1}) + A(\Delta u_iv_iu_{i+1})
+ A (\Delta u_{i+1}v_{i+1}u_i) +A(\Delta u_{i+1}v_{i+1}v_i))\\
&\leq& \frac{1}{2}(|v_iv_{i+1}|+ |u_iu_{i+1}|) \frac{\sqrt{2}}{2} \\
& \leq & \frac{\sqrt{2}}{4}(|v_iv_{i+1}|+|u_iv_i|+|v_iv_{i+1}|+|v_{i+1}u_{i+1}|)\\
 &\leq& \frac{\sqrt{2}}{4}(2|v_iv_{i+1}|+ \sqrt{2})\\
&=&  \frac{\sqrt{2}}{2}|v_iv_{i+1}| +\frac{1}{2}.
\end{eqnarray*}
Summing up, we get
\begin{eqnarray*}
  |A(Q)-A(v+tP)| &\leq& \sum_{i=1}^m A(S_i) \\
 &\leq& \sum_{i=1}^m(\frac{\sqrt{2}}{2}|v_iv_{i+1}| +\frac{1}{2})\\
&=& \frac{\sqrt{2}}{2}Lt + \frac{m}{2}.
\end{eqnarray*}

Let $T_i$ be the set of grid points inside $S_i$ or on the line
segment $u_iu_{i+1}$. Let $P_i$ be the convex set spanned by $T_i$.
Applying Pick's theorem  to $P_i$, we have
$$A(P_i)=I(P_i)+\frac{B(P_i)}{2}-1.$$
Thus
\begin{eqnarray*}
  |T_i| &=& I(P_i) + B(P_i)\\
     &\leq & 2(A(P_i)+1)\\
     &\leq& 2(A(S_i)+1)\\
   &\leq& \sqrt{2}|v_iv_{i+1}| +3.
\end{eqnarray*}
Summing up, we get
\begin{eqnarray*}
  |I(Q)-I(v+tP)|&\leq& \sum_{i=1}^m |T_i|\\
&\leq& \sum_{i=1}^m \sqrt{2}|v_iv_{i+1}| +3 \\
&=& \sqrt{2}tL+3m.
\end{eqnarray*}
Putting together, we have
\begin{eqnarray*}
  |I(v+tP)-A(v+tP)| \hspace*{-2cm}&&\\
&\leq& | I(Q)-A(Q)| + |A(Q)-A(v+tP)| + |I(Q)-I(v+tp)|\\
&\leq&   \frac{1}{2}(tL+(\sqrt{2}+1)m)-1 + (\frac{\sqrt{2}}{2}tL + \frac{m}{2})
+ (\sqrt{2}tL+3m)\\
&=& \frac{3\sqrt{2}+1}{2}Lt+ (4+\frac{\sqrt{2}}{2})m-1\\
&<& 3Lt +5m.
\end{eqnarray*}
The proof of this lemma is finished. \hfill $\square$


By counting grid points in $(-x_0/b, -y_0/b)+ (n/b)P$, we get the
following corollary.  Since the number of grid points on the
boundary of $nP$ is always a lower order term,  it does not matter
whether grid points on the boundary are included or not. In fact, in
 latter applications, the polygon region $P$ is often defined with
one part of  boundary included while the other part of boundary excluded.

\begin{corollary}\label{lattice}
For any fixed point $(x_0,y_0)$, let $L_b$ be a lattice $\{(x_0+ib,y_0+jb)\colon i,j
\in \Z\}$ and $P$ be a simple polygon. For $n\gg b$, we have
$$|L_b\cap nP|= \frac{n^2}{b^2}A(P)+O(\frac{n}{b}).$$
\end{corollary}

\noindent
{\bf Proof of Lemma \ref{l0}:} Let $B$ be a ``good''
$2$-coloring/bit-string of $\Z_b$ with $m_k(\Z_b)b^2$ monochromatic
$k$-APs. Any $k$-AP in $\Z_b$ can be parametrized by a pair
$(a',d')$ satisfying $0\leq a', d'\leq b-1$. Let $S$ be the set of
parameters $(a',d')$ such that the corresponding $k$-APs in $\Z_b$
are monochromatic. We have
$$|S|=m_k(\Z_b)b^2.$$

For sufficiently large $n$, we write $n=bt+r$ with $0\leq r\leq
b-1$. Consider the periodic construction $BB\cdots BR$ (see
(\ref{eq:pcon})). Note that the number of $k$-APs containing some
elements of $R$ is $O(n)$. We need  estimate monochromatic $k$-APs
lying entirely in  $[bt]$.

Let $P$ be a parallelogram defined by
$$P=\{(x,y)\colon 0<  x \leq 1 \mbox{ and } 0< x+y(k-1)\leq 1\}.$$
The area of $P$ is clearly $\frac{1}{(k-1)}$.

A $k$-AP $(a, a+d,  \ldots, a+(k-1)d)$  in $[bt]$ is monochromatic if and only
if $(a \mod b, d \mod b)\in S$.
Let  $L_b^{(a',d')}$ be the lattice $\{(a'+ib, d'+jb)\colon i,j\in \Z\}$.
Applying Corollary \ref{lattice},
 the number of monochromatic $k$-APs in $[bt]$
is exactly
\begin{eqnarray*}
  \sum_{(a',d')\in S} |L_b^{(a',d')}\cap (bt)P| &=& \sum_{(a',d')\in S} A(P)t^2 +O(t)\\
&=& |S|A(P)t^2 +O(b^2t)\\
&=&\frac{1}{k-1} m_k(\Z_b)(bt)^2+O(b^2t).
\end{eqnarray*}
Thus,
$$m_k([n],c)=\frac{1}{k-1} m_k(\Z_b)(bt)^2+O(b^2t).$$
Note that the  number of $k$-APs in $[n]$ is $\frac{n^2}{k-1}+O(n)$.
Taking the ratio, we have
$$m_k([n])\leq \frac{m_k([n],c)}{AP([n])}=m_k(\Z_b) +O(\frac{1}{t}).$$
First taking (upper) limit as $n$ goes to infinity, we get
 $$\varlimsup_{n\to\infty} m_k([n])\leq m_k(\Z_b) .$$
Then taking (lower) limit as $b$ goes to infinity, we have
$$\varlimsup_{n\to\infty} m_k([n])\leq \varliminf_{b\to\infty}m_k(\Z_b). $$
The proof of Lemma \ref{l0} is finished. 
\hfill $\square$

\subsection{Proof of Theorem  \ref{t:m4ub}}


It suffices to consider  the case that $n$ is not divisible by $b$. Write
$n=bt+r$ with $1\leq r\leq b-1$.  Recall the periodic
construction
$$\underbrace{BB\cdots B}_tR.$$
Here $R$ is any bit-string of length $r$.

The number of  $4$-APs containing some bit(s) in $R$ is $O(n)$. The
major term in the number of all monochromatic $4$-APs depends only
on $B$ and $r$. We divide the set of all non-degenerated $4$-APs in
$\Z_n$ into eight classes $C_i$ for $0\leq i \leq 7$.
\begin{center}
\begin{tabular}{|c |c |}\hline
Classes  & the meaning in $\mathbb Z_n$ \\
\hline
$C_0$ &  $ a<a+d<a+2d<a+3d<n$ \\
\hline
$C_1$ & $a<a+d<a+2d<n \leq a+3d<2n$
\\ \hline
$C_2$ & $a<a+d<n \leq  a+2d<a+3d<2n$
\\ \hline
$C_3$ & $ a<a+d<n \leq a+2d<2n \leq a+3d<3n$
\\ \hline
$C_4$ &  $ a<n \leq a+d<a+2d<a+3d<2n$
\\ \hline
$C_5$ & $a<n \leq a+d<a+2d<2n \leq a+3d<3n$   \\
\hline
$C_6$ & $a<n \leq  a+d<2n \leq a+2d< a+3d<3n$
\\ \hline
$C_7$ & $ a<n \leq a+d<2n \leq a+2d<3n \leq a+3d<4n$
\\ \hline
\end{tabular}
\end{center}

\begin{figure}[htbp]
 \centerline{ \psfig{figure=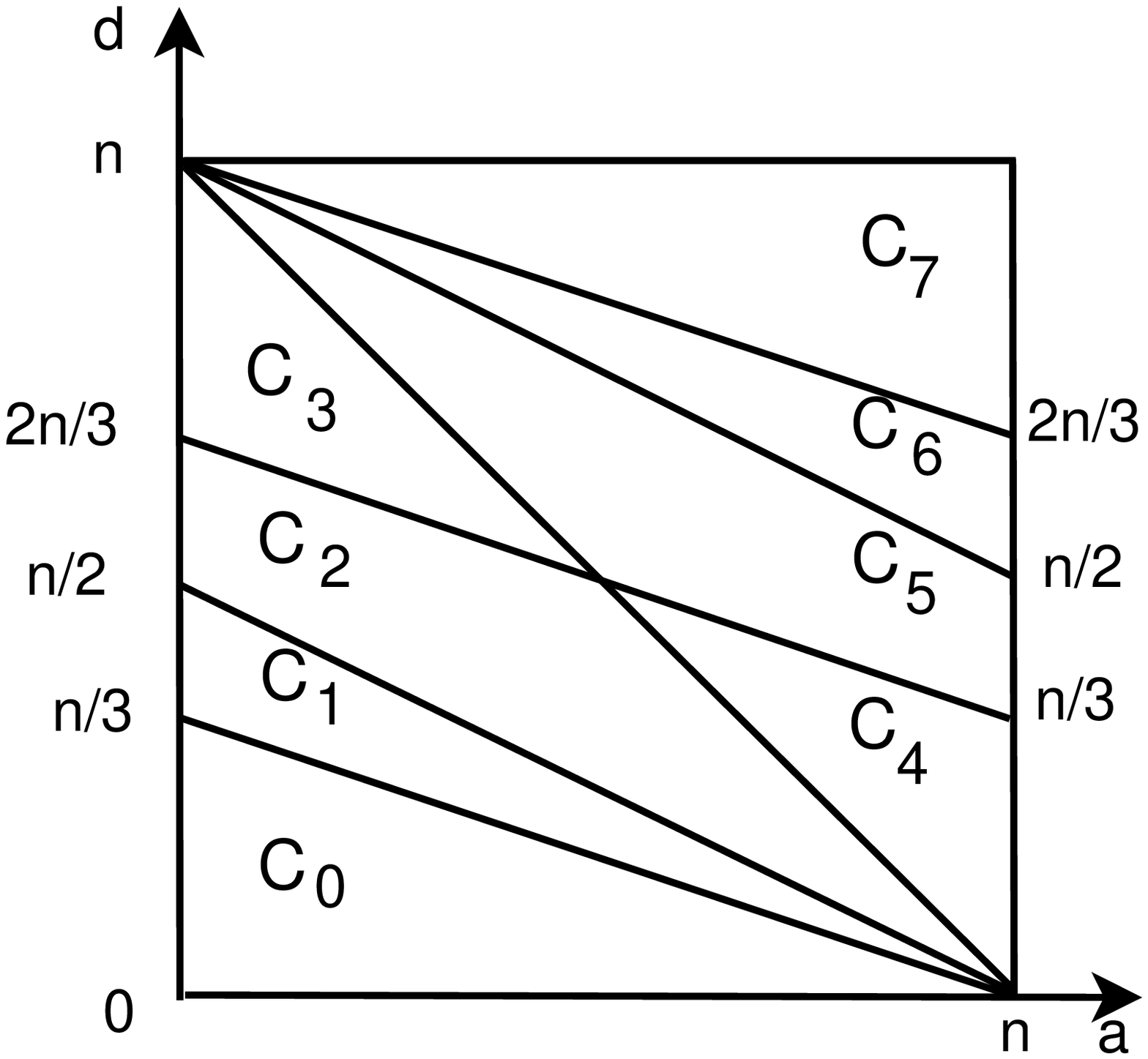, width=0.4\textwidth} \hfil  \psfig{figure=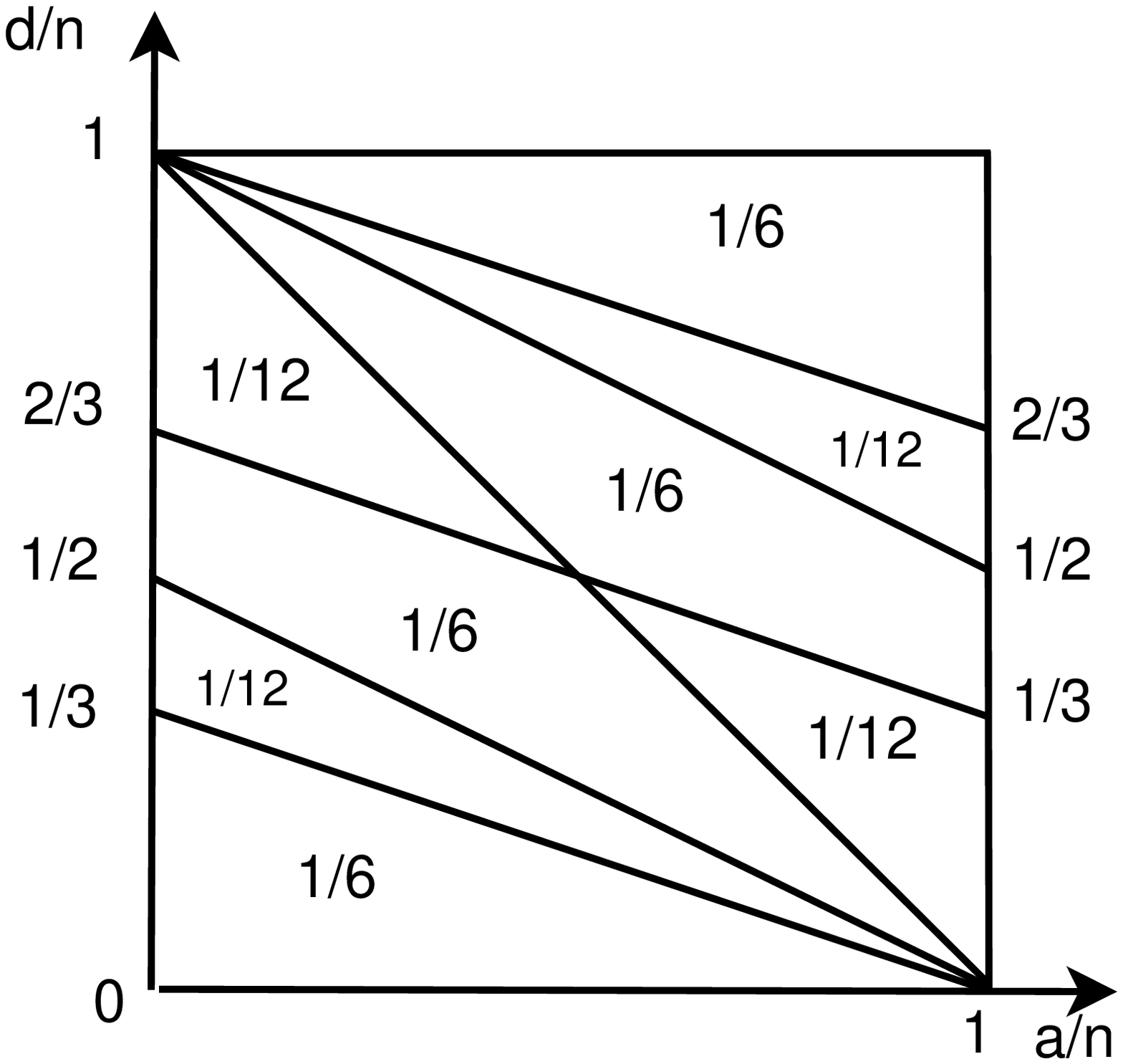, width=0.4\textwidth} }
\begin{multicols}{2}
\caption{{\it The eight regions in the parameter space of all 4-APs in
$\mathbb Z_n$. }
  \label{regions}} \newpage
\caption{{\it The areas of the eight  normalized regions.}
 \label{areas}}
\end{multicols}
\end{figure}

 These 8 classes can be viewed as 8 regions in the parameter space of $(a,d)$
as shown in Figure \ref{regions}. Let us normalize the parameters so
that the area of the whole square is $1$.
 For $0\leq i\leq 7$, let $a_i$ be the area of the $i$-th normalized region. We have
$$a_0=\frac{1}{6}, \; a_1=\frac{1}{12}, \; a_2=\frac{1}{6}, \; a_3=\frac{1}{12},  \;
a_4=\frac{1}{12}, \; a_5=\frac{1}{6}, \; a_6=\frac{1}{12},  \; a_7=\frac{1}{6}.$$

For $r_1, r_2, r_3\geq 0$, an {\it $(r_1, r_2, r_3)$-generalized 4-term arithmetic progression}
is of form
$$a, a+d-r_1,a+2d-(r_1 +r_2),a+3d-(r_1 +r_2+r_3).$$
Here $(a,d)$ are the parameters determining the
$(r_1, r_2, r_3)$-generalized 4-term arithmetic progression.

We have the following Lemma.
\begin{lemma} \label{upper}
  For $0\leq i\leq 7$, write $i$ as a bit-string $x_1x_2x_3$ of length three.
Let $c_i$ be the number of all monochromatic    $(x_1r, x_2r, x_3r)$-generalized 4-term arithmetic progressions in $B$.
Then the number of monochromatic $4$-APs  in  $BB\cdots BR$
is $$\sum_{i=0}^7a_ic_it^2 +O(t).$$
In particular, we have
$$m_4(\Z_n)\leq \sum_{i=0}^7a_i\frac{c_i}{b^2}+o(1).$$
\end{lemma}
{\bf Proof:} A $4$-AP is said on {\it the boundary of some $C_i$} if it is in $C_i$ and
 contains an element in $R$. Note that the number of $4$-APs on the boundary is $O(n)$. We can ignore
these $4$-APs in the calculation below.

For any $0\leq a', d'\leq b-1$,
the lattice $L^{(a',d')}_b =\{(a'+ub, d'+vb)\colon 0\leq u, v<t\}$ distributes
evenly in the square $[0,n)\times [0,n)$. Applying corollary \ref{lattice},
 we have
$$|L^{(a',d')}_b\cap C_i| =a_i t^2 +O(t)$$
for $0 \leq i \leq 7.$ We also observe any monochromatic 4-term
arithmetic progression with parameter $(a,d)=(a'+ub, d'+vb)\in C_i$
if and only if the $(x_1r, x_2r, x_3r)$-generalized 4-term
arithmetic progression with parameter $(a',d')$ is monochromatic in
$B$. Thus the number of
 monochromatic  4-term arithmetic progressions with parameter $(a,d)\in C_i$
is $$c_ia_it^2+O(t).$$ Hence the number of  monochromatic  4-term
arithmetic progressions in $BB\cdots BR$ is
$\sum_{i=0}^7a_ic_it^2+O(t)$ and $m_4(\Z_n)\leq
\sum_{i=0}^7a_ic_i/b_i^2 +O(1/n)$. \hfill $\square$

We are ready to prove Theorem \ref{t:m4ub}.

\noindent
{\bf Proof of Theorem \ref{t:m4ub}:} Recall $B_{20}=(1,1,1,0,1,1,0,1,1,1,0,0,0,1,0,0,1,0,0,0)$ and $B_{22}=(1,1,1,0,1,1,0,1,0,0,0,1,1,1,0,1,0,0,1,0,0,0).$

When $n$ is odd, we use the periodic construction $B_{20}B_{20}\cdots
B_{20}R$.
 We write $n=20t+r$, where $r=1,3,5,\ldots, 19$. For each odd $r$,  it turns out that
the values of  $c_i$ depends only on $i$ but not on $r$. These
values are  given in Table \ref{tab:3}.

\begin{table}[htbp]
\begin{center}
\begin{tabular} {|c|c|c|c|c|c|c|c|}\hline
   $c_0$ & $c_1$ & $c_2$ & $c_3$ & $c_4$  & $c_5$  & $c_6$  & $c_7$\\ \hline
  36    & 50    & 50    & 50    & 50     & 50     & 50 &36\\
\hline
\end{tabular}
\caption{The values of $c_i$'s for $B_{20}$ and  any odd $r$  satisfying $1 \leq r
\leq 19$.} \label{tab:3}
\end{center}
\end{table}
By Lemma \ref{upper}, we have
$$m_{4}(\Z_n)\leq \sum_{i=0}^7\frac{a_ic_i}{b^2}+o(1)
=\frac{17}{150}+o(1).$$

When $n$ is even,  we prove  only a weaker result $m_4(\Z_n)\leq
\frac{175}{1452}+o(1)< 0.12052342$ here and {\em postpone the proof
of actual bound until the end of next section}. We write $n=22t+r$
where $r=0,2,4,\ldots, 20$. We use the periodic construction
$B_{22}B_{22}\cdots B_{22}R$. If $r=0$, then we have
$$m_4(\Z_n)\leq m_4(\Z_{22})=\frac{21}{242}<0.086777.$$
We are done in this case. For $r=2,4,6,\ldots, 20$, the values
$c_i$ depends only on $i$, but not on $r$. These values are  given
in Table \ref{tab:4}.

\begin{table}[htbp]
\begin{center}
\begin{tabular} {|c|c|c|c|c|c|c|c|}\hline
 $c_0$ & $c_1$ & $c_2$ & $c_3$ & $c_4$  & $c_5$  & $c_6$  & $c_7$\\ \hline
 42    & 63    & 70    & 63    & 63     & 70     & 63 &42\\
\hline
\end{tabular}
\caption{The values of $c_i$'s for $B_{22}$ and  each even $r$ such that $2 \leq r
\leq 20$.} \label{tab:4}
\end{center}
\end{table}
By Lemma \ref{upper}, we have
$$m_{4}(\Z_n)\leq \sum_{i=0}^7\frac{a_ic_i}{b^2}+o(1)
=\frac{175}{1452}+o(1). $$
\hfill $\square$

\subsection{The Proof of Theorem \ref{t:m5ub} }


Let $B$ be a "good" 2-coloring of $\Z_b$.  We consider the periodic
construction $c=BB\cdots BR$.

\begin{figure}[htbp]
 \centerline{ \psfig{figure=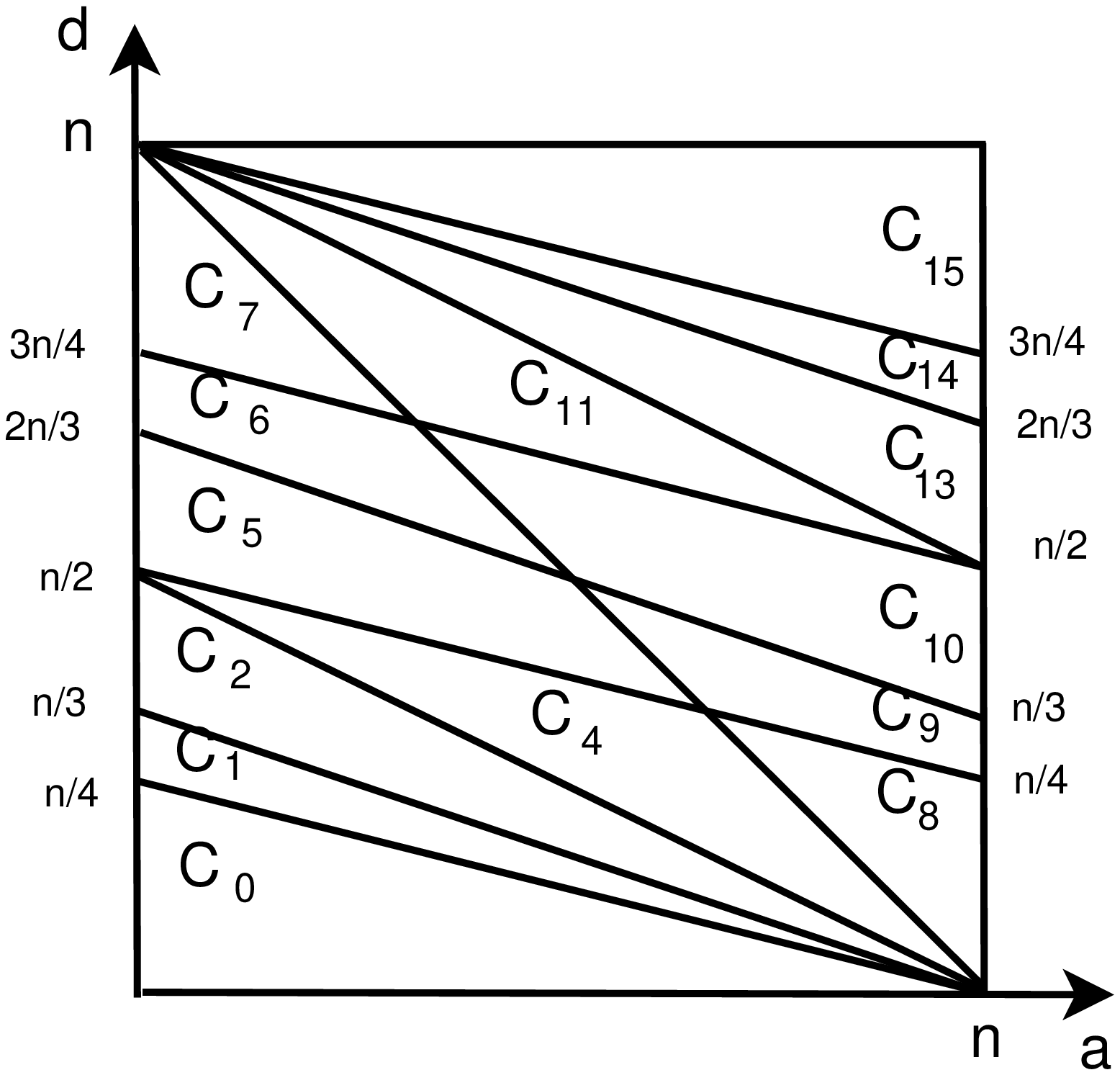, width=0.4\textwidth} \hfil  \psfig{figure=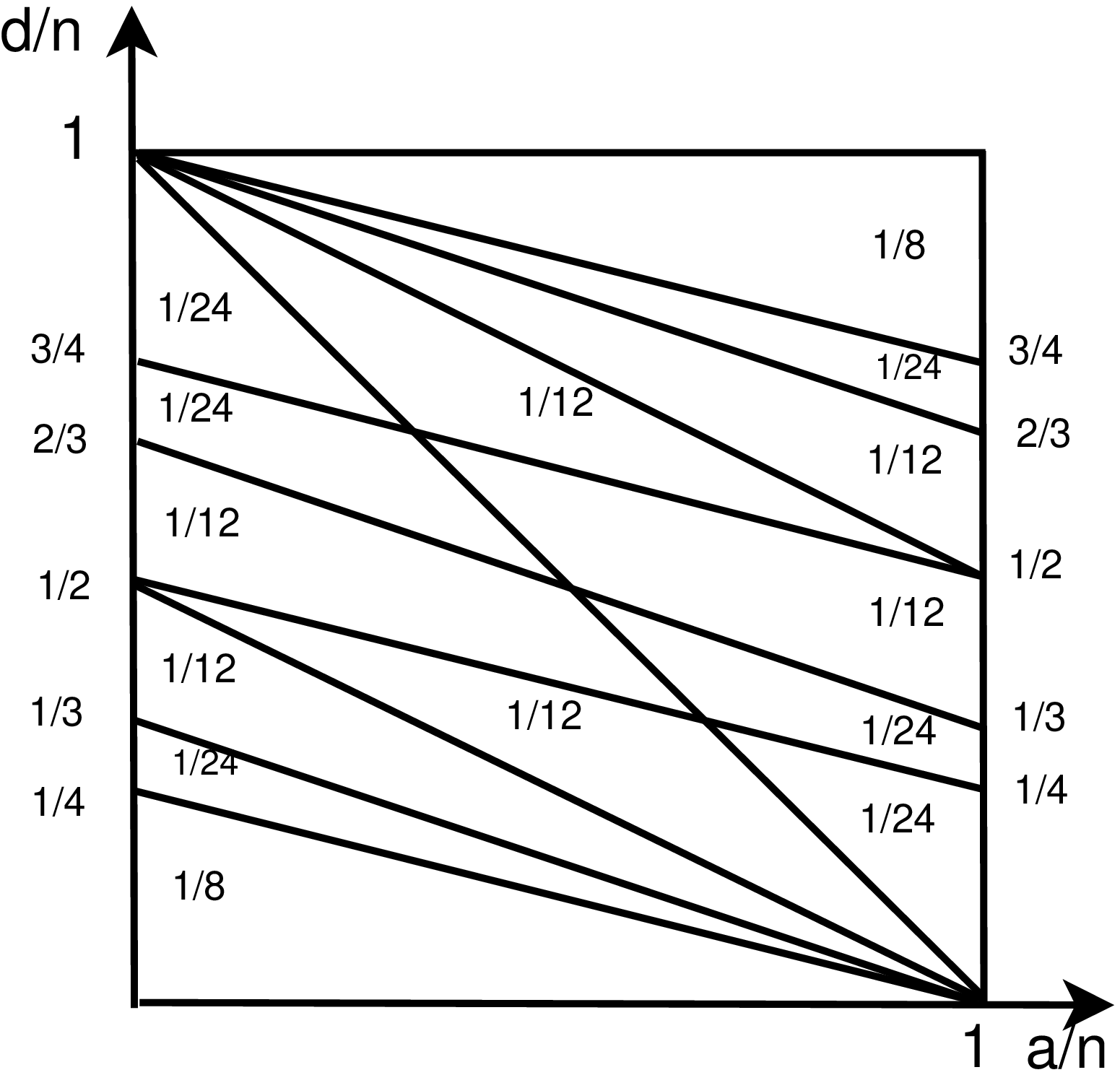, width=0.4\textwidth} }
\begin{multicols}{2}
\caption{{\it The 14 regions of the parameter space of all 5-APs in $\mathbb Z_n$. }
  \label{regions:5}} \newpage
\caption{{\it The areas of the 14 normalized regions.}
 \label{areas:5}}
\end{multicols}
\end{figure}

Similar to the proof of Theorem \ref{t:m4ub}, we can divide all
non-degenerated $5$-APs into 14 classes $C_i$ with index $i$ in $S
=\{0,\ldots,15\} \setminus \{3,12\}$, see table below and Figure
\ref{regions:5}; let $a_i$ be the area of $i$-th normalized region
(see Figure \ref{areas:5}). We have
$$a_0=\frac{1}{8}, \; a_1=\frac{1}{24},\; a_2=\frac{1}{12}, \; a_4=\frac{1}{12},\; 
a_5=\frac{1}{12},\; a_6=\frac{1}{24}, \; a_7=\frac{1}{24},$$ 
$$a_8=\frac{1}{24}, \; a_9=\frac{1}{24},\; a_{10}=\frac{1}{12}, \; a_{11}=\frac{1}{12},\; 
a_{13}=\frac{1}{12},\; a_{14}=\frac{1}{24}, \; a_{15}=\frac{1}{8}.$$ 
\begin{center}
\begin{tabular}{|c |c |}\hline
type  & the meaning in $\mathbb Z_p$ \\
\hline

$C_0$ &  $ a<a+d<a+2d<a+3d<a+4d<n$ \\
\hline

$C_1$ & $a<a+d<a+2d<a+3d<n \leq a+4d<2n$
\\ \hline

$C_2$ & $a<a+d< a+2d<n \leq a+3d<a+4d<2n$
\\ \hline

$C_4$ &  $ a< a+d<n \leq a+2d<a+3d<a+4d<2n$
\\ \hline

$C_5$ & $a< a+d<n
\leq a+2 d < a+3d < 2 n \leq a+4d<3n$   \\
\hline

$C_6$ & $a< a+d<n  \leq a+2d<2n \leq  a+3d<a+4d<3n$
\\ \hline

$C_7$ & $ a<a+d<n \leq a+2d<2n \leq a+3d<3n \leq a+4d<4n$
\\ \hline

$C_8$ & $ a<n \leq a+d<a+2d<a+3d<a+4d<2n$
\\ \hline

$C_9$ & $ a<n  \leq a+d<a+2d<a+3d<2n \leq a+4d<3n$
\\ \hline

$C_{10}$ & $ a<n \leq a+d<a+2d <2n \leq a+3d<a+4d<3n$
\\ \hline
$C_{11}$ & $ a<n \leq a+d<a+2d<2n \leq a+3d<3n \leq a+4d<4n$
\\ \hline
$C_{13}$ & $ a< n \leq a+d<2n \leq a+2d<a+3d<3n \leq a+4d<4n$
\\ \hline
$C_{14}$ & $ a<n  \leq a+d<2n \leq a+2d < 3n \leq a+3d<a+4d<4n$
\\ \hline
$C_{15}$ & $ a<n \leq a+d<2n \leq a+2d<3n \leq a+3d<4n \leq a+4d<5n$
\\ \hline
\end{tabular}
\end{center}

Assume $r_i \geq 0$ for $1 \leq i \leq 4$. An {\it
$(r_1,r_2,r_3,r_4)$-generalized 5-term arithmetic progression} is of
form
$$
a,a+d-r_1,a+2d-(r_1+r_2),a+3d-\sum_{i=1}^3 r_i,a+4d-\sum_{i=1}^4
r_i.
$$
Given $(r_1,r_2,r_3,r_4)$, an $(r_1,r_2,r_3,r_4)$-generalized 5-term
arithmetic progression is determined by $(a,d)$. We have the
following lemma; we will omit the proof since it is similar to Lemma
\ref{upper}.

\begin{lemma}\label{le3}
Let  $S=\{0,\ldots,15\} \setminus
\{3,12\}$. For each $i \in S$, write $i$ as a bit-string
$x_1x_2x_3x_4$ of length four. Let $c_i$ be the number of all
monochromatic $(x_1r, x_2r, x_3r,x_4r)$-generalized 5-term
arithmetic progressions in $B$. Then the number of monochromatic
$5$-APs  in  $BB\cdots BR$ is
$$\sum_{i \in S} a_ic_it^2 +O(t).$$ In particular, we have
$$m_5(\Z_n)\leq \sum_{i \in S} a_i\frac{c_i}{b^2}+o(1).$$
\end{lemma}


\noindent
{\bf Proof of Theorem \ref{t:m5ub}:} Recall  the 2-coloring $B_{74}$
of $\mathbb Z_{74}$ as following

\noindent$\{
1,1,1,1,0,1,1,1,0,0,0,0,1,0,1,1,0,0,1,0,1,0,1,0,0,1,1,0,1,0,0,0,0,1,1,1,0,\\
 \hspace*{0.2cm}  1,1,1,1,0,1,1,1,0,0,0,0,1,0,1,1,0,0,1,0,0,0,1,0,0,1,1,0,1,0,0,0,0,1,1,1,0\}.$

 Write $n=74t+r$, where $0 \leq r \leq 73$.
If $r=0$, then $m_5(\Z_n)\leq m_5(\Z_{74})=\frac{73}{2738}$; we are
done in this case. Now we assume $r\not=0$ and use the periodic
construction $B_{74}B_{74}\cdots B_{74}R$, where $R$ is any
bit-string of length $r$.

The values $c_i$ in Lemma
\ref{le3} depend on $i$ and $r$.  These values are given
in Table \ref{tab:5}.

\begin{table}[htbp]
\begin{center}
\small
\begin{tabular}{|c|c|c|c|c|c|c|c|c|c|c|c|c|c|c|c|} \hline
\!\!
values \!\!\!\!&\!\!\!\! $c_0$ \!\!\!\!&\!\!\!\! $c_1$ \!\!\!\!&\!\!\!\! $c_2$ \!\!\!\!&\!\!\!\! $c_4$ \!\!\!\!&\!\!\!\! $c_5$    \!\!\!\!&\!\!\!\! $c_6$   \!\!\!\!&\!\!\!\! $c_7$ \!\!\!\!&\!\!\!\! $c_8$ \!\!\!\!&\!\!\!\! $c_9$  \!\!\!\!&\!\!\!\! $c_{10}$     \!\!\!\!&\!\!\!\! $c_{11}$    \!\!\!\!&\!\!\!\! $c_{13}$    \!\!\!\!&\!\!\!\! $c_{14}$ \!\!\!\!&\!\!\!\! $c_{15}$ \\ \hline
\!\!
even $r\not=0$ \!\!\!\!&\!\!\!\!146 \!\!\!\!&\!\!\!\! 293   \!\!\!\!&\!\!\!\! 377   \!\!\!\!&\!\!\!\! 377   \!\!\!\!&\!\!\!\! 378      \!\!\!\!&\!\!\!\! 359   \!\!\!\!&\!\!\!\!  293   \!\!\!\!&\!\!\!\!  293   \!\!\!\!&\!\!\!\! 359     \!\!\!\!&\!\!\!\! 378          \!\!\!\!&\!\!\!\! 377         \!\!\!\!&\!\!\!\! 377         \!\!\!\!&\!\!\!\! 293       \!\!\!\!&\!\!\!\! 146\\
\hline
\!\!
odd $r\not=37$ \!\!\!\!&\!\!\!\!
146  \!\!\!\!&\!\!\!\!  293   \!\!\!\!&\!\!\!\! 375   \!\!\!\!&\!\!\!\! 375   \!\!\!\!&\!\!\!\! 374      \!\!\!\!&\!\!\!\! 357     \!\!\!\!&\!\!\!\!  293  \!\!\!\!&\!\!\!\!293   \!\!\!\!&\!\!\!\! 357    \!\!\!\!&\!\!\!\! 374         \!\!\!\!&\!\!\!\!  375        \!\!\!\!&\!\!\!\!  375        \!\!\!\!&\!\!\!\! 293     \!\!\!\!&\!\!\!\! 146\\ \hline
\!\! $r=37$ \!\!\!\!&\!\!\!\!146 \!\!\!\!&\!\!\!\!144  \!\!\!\!&\!\!\!\!144\!\!\!\!&\!\!\!\! 144\!\!\!\!&\!\!\!\! 144 \!\!\!\!&\!\!\!\!144 \!\!\!\!&\!\!\!\!144\!\!\!\!&\!\!\!\!144 \!\!\!\!&\!\!\!\!144\!\!\!\!&\!\!\!\! 144\!\!\!\!&\!\!\!\! 144 \!\!\!\!&\!\!\!\!144\!\!\!\!&\!\!\!\! 144\!\!\!\!&\!\!\!\! 146 \\ \hline
\end{tabular}
\caption{The values of $c_i$'s for $B_{74}$ and various $r$'s.}
\label{tab:5}
\end{center}
\end{table}
By Lemma \ref{le3}, we have
$$m_5(\Z_n) \leq \sum_{i \in S} a_i \frac{c_i}{b^2}+o(1)=
\left\{
  \begin{array}[c]{ll}
    3629/65712 +o(1) & \mbox{ for odd } r \not=37,\\
   289/10952 +o(1) &  r =37,\\
   3647/65712 +o(1) & \mbox{ for even } r \not=0.\\
  \end{array}
\right.
$$
The proof of Theorem \ref{t:m5ub} is finished. \hfill $\square$

\section{Proof of Theorem \ref{t:m5ub}} 
Recall the good coloring of $\Z_{22}$:
\[B_{22}=(1,1,1,0,1,1,0,1,0,0,0,1,1,1,0,1,0,0,1,0,0,0).\]
We observe that the first 11 coordinates and last 11 coordinates
differ only by $1$ bit. Let $B_{11}=(1,1,1,0,1,\ast,0,1,0,0,0)$,
where `$\ast$' could be either $0$ or $1$. $B_{11}$ has the
following property. (This is because  $B_{22}$
contains no non-degenerated monochromatic $4$-APs.)
\begin{property} \label{b11}
No matter which bit-value the  `$\ast$' takes,
$B_{11}$ contains no non-degenerated monochromatic
$4$-APs of  $\Z_{11}$.
\end{property}

\begin{lemma}\label{recursive11}
For any $t\geq 2$, we have
 \begin{equation}
   \label{eq:recursive11}
   m_4(\Z_{11t})\leq \frac{10+m_4(\Z_{t})}{121}.
 \end{equation}
\end{lemma}
{\bf Proof:} Let $B_t$ be a 2-coloring/bit-string of $\Z_t$ which
has exactly  $m_4(\Z_t)t^2$ monochromatic $4$-APs.
First we consider the periodic construction
 $$\underbrace{B_{11}B_{11}\cdots B_{11}}_t.$$
Each block $B_{11}$ has exactly one `$\ast$'; there are $t$ $\ast$'s
in total. Finally, we replace these $\ast$'s by the values of $B_t$
in the cyclic order. We denote the coloring by $B_{11}\ltimes B_t$.
(For example, $B_{22}=B_{11}\ltimes (1,0)$.)

Because of Property \ref{b11},  a $4$-AP of $B_{11}\ltimes B_t$ with parameter $(a,d)$
is monochromatic only if $11\mid d$. The number of monochromatic  $4$-APs
is exactly $10t^2+m_4(\Z_t)t^2.$ We have
$$m_4(\Z_{11t}) \leq \frac{10t^2+m_4(\Z_t)t^2}{(11t)^2}= \frac{10+m_4(\Z_t)}{121}.$$
The proof of this lemma is finished.  \hfill $\square$

A similar construction can be applied to $5$-APs. Let \\
{\small $B_{37}=( 1,1,1,1,0,1,1,1,0,0,0,0,1,0,1,1,0,0,1,0,\ast,0,1,0,0,1,1,0,1,0,0,0,0,1,1,1,0).$}

Let $B_t$ be a 2-coloring/bit-string of $\Z_t$ which
has exactly  $m_5(\Z_t)t^2$ monochromatic $4$-APs.
 We can define $B_{37}\ltimes B_t$ similarly.
For example, $B_{74}=B_{37}\ltimes (1,0)$.
Note that $B_{74}$ contains no non-degenerated monochromatic
$5$-APs. We have the following property. 
\begin{property} \label{b37}
No matter which bit-value the `$\ast$' takes,
$B_{37}$ contains no non-degenerated monochromatic
$5$-APs of  $\Z_{37}$.
\end{property}
Using this property and the construction $B_{37}\ltimes B_t$,
we have the following lemma. The proof is omitted.

\begin{lemma} \label{recursive37}
For any $t\geq 2$, we have
 \begin{equation}
   \label{eq:recursive37}
   m_5(\Z_{37t})\leq \frac{36+m_4(\Z_{t})}{37^2}.
 \end{equation}
\end{lemma}

\noindent
{\bf Proof of Theorem \ref{tzb}:}
Applying Lemma \ref{recursive11} recursively, we have
\begin{eqnarray*}
  m_4(\Z_{11^{s}}) &\leq& \frac{10}{11^2} + \frac{1}{11^2}m_4(\Z_{11^{s-1}}) \\
&\leq& \frac{10}{11^2} + \frac{10}{11^4} + \frac{1}{11^4}m_4(\Z_{11^{s-2}})\\
&\leq& \cdots\\
&\leq&  \frac{10}{11^2} + \frac{10}{11^4} + \cdots + \frac{10}{11^{2s}} +
\frac{1}{11^{2s}}m_4(\Z_{1})\\
&=& \frac{10}{11^2} \frac{1-\frac{1}{11^{2s}}}{1-\frac{1}{11^2}}+\frac{1}{11^{2s}}\\
&=& \frac{1}{12} + \frac{1}{12\times 11^{2s-1}}.
\end{eqnarray*}
Thus,
\begin{eqnarray*}
\varliminf_{n\to \infty} m_4(\Z_n)&\leq&
 \varliminf_{s\to\infty}  m_4(\Z_{11^{s}})\\
&\leq& \lim_{s\to\infty} (\frac{1}{12} + \frac{1}{12\times 11^{2s-1}})\\
&=& \frac{1}{12}.
\end{eqnarray*}
Similarly, from Lemma \ref{recursive37}, we can show
$\varliminf_{n\to \infty} m_5(\Z_n)\leq \frac{1}{38}.$
\hfill $\square$

\noindent {\bf Proof  of Theorem \ref{t:m4ub} for even $n$:} Here we
assume $n$ is even and not divisible by $22$. Let $B=B_{11}\ltimes
B_{20}$ which is a $2$-coloring of $\Z_{220}$. Write $n=220t+r$ with
$0\leq r\leq 218$. Here $r$ is even and not divisible by $22$.
Consider the periodic construction $BB\cdots BR$ as before. For
these $r$, the values  $c_i$ depends only on $i$ but not on $r$.
These values are given in Table \ref{tab:6}.

\begin{table}[htbp]
\begin{center}
\begin{tabular} {|c|c|c|c|c|c|c|c|}\hline
 $c_0$ & $c_1$ & $c_2$ & $c_3$ & $c_4$  & $c_5$  & $c_6$  & $c_7$\\ \hline
 4882 &7563 &8230 &7563 &7563 &8230 &7563& 4882\\
\hline
\end{tabular}
\caption{The values of $c_i$'s for $B=B_{11}\ltimes B_{20}$ and
 even $r=2,4,\ldots, 218$ such that $r$ is not divisible by $22$.} \label{tab:6}
\end{center}
\end{table}
By Lemma \ref{upper}, we have
$$m_{4}(\Z_n)\leq \sum_{i=0}^7\frac{a_ic_i}{b^2} +o(1)
=\frac{8543}{72600}+o(1)< 0.11767722. $$ 
The remaining case is proved.\hfill $\square$

\section{Proof of  Theorem \ref{t:m4lb}}
In this section, we will deal with lower bound of $m_4(\Z_n)$.

\noindent
 {\bf Proof of Theorem \ref{t:m4lb}:} Given a 2-coloring $c$ of $\mathbb Z_n$, we will establish
 an
inequality which is similar to equation (4.8) in \cite{cameron}.
For each $0 \leq i \leq 4$, let $u_i$ be the number of $4$-APs with
exactly $i$ red numbers.  We have
\begin{eqnarray*}
u_1+u_3& =& |A_1\cap B_2 \cap B_3 \cap B_4|+|B_1\cap A_2 \cap B_3
\cap B_4|+|B_1\cap B_2 \cap A_3 \cap B_4|\\
&+&|B_1\cap B_2 \cap B_3 \cap A_4|+ |B_1\cap A_2 \cap A_3 \cap
A_4|+|A_1\cap B_2 \cap A_3 \cap A_4|\\
&+&|A_1\cap A_2 \cap B_3 \cap A_4|+|A_1\cap A_2 \cap A_3 \cap B_4|.
\end{eqnarray*}
Note that
\begin{equation}
  \label{eq:ie}
|A_1\cap B_2 \cap B_3 \cap B_4|=|B_2 \cap B_3 \cap B_4|-|B_1 \cap B_2 \cap B_3 \cap B_4|.
\end{equation}
Applying equations similar to  (\ref{eq:ie}), we get
\begin{equation}
\label{eq:ab0}
4u_0+ u_1+u_3 +4u_4 = \sum_{1 \leq i < j <k \leq 4} (|A_i \cap A_j \cap A_k|+|B_i \cap
B_j \cap B_k|).
\end{equation}

By the inclusion-exclusion formula,  for any $1 \leq i < j <k \leq 4$, we
have
$$
|A_i \cup A_j \cup A_k|= \sum_{s \in \{i,j,k\}}|A_s|-\sum_{\{s,t\}
\in \binom{\{i,j,k\}}{2}} |A_s \cap A_t|+|A_i \cap A_j \cap A_k|.
$$
Since $|A_i \cup A_j \cup A_k|=n^2-|B_i \cap B_j \cap B_k|$, we have
\begin{equation}
  \label{eq:ab1}
  |A_i \cap A_j \cap A_k|+|B_i \cap
B_j \cap B_k|=n^2 - \sum_{s \in \{i,j,k\}}|A_s|+\sum_{\{s,t\} \in
\binom{\{i,j,k\}}{2}} |A_s \cap A_t|.
\end{equation}
By the symmetry of $A_i$'s and $B_i$'s, we get
\begin{equation}
  \label{eq:ab2}
  |A_i \cap A_j \cap A_k|+|B_i \cap
B_j \cap B_k|=n^2 - \sum_{s \in \{i,j,k\}}|B_s|+\sum_{\{s,t\} \in
\binom{\{i,j,k\}}{2}} |B_s \cap B_t|.
\end{equation}
Combining equations (\ref{eq:ab1}) and (\ref{eq:ab2}) and summing
over $1\leq i<j<k\leq 4$, we get
\begin{eqnarray*}
 2\!\!\!\! \sum_{1 \leq i < j <k \leq 4} ( |A_i \cap A_j \cap A_k|+|B_i \cap
B_j \cap B_k|) \hspace*{-6cm} &&\\
&=&  \sum_{1 \leq i < j <k \leq 4} \left (2n^2 -\!\!\!\! \sum_{s \in
\{i,j,k\}}\!\!\!\!(|A_s|+ |B_s|)+ \!\!\!\!\!\!\!\! \sum_{\{s,t\}
\in \binom{\{i,j,k\}}{2}} \!\!\!\!\!\!\!\! (|A_s\cap A_t|+ |B_s \cap B_t|) \right)\\
&=& 8n^2- 12n^2 + 2\sum_{1\leq i < j \leq 4 }(|A_i\cap A_j|+|B_i \cap
B_j|)\\
&=& -4n^2+ 2\sum_{1\leq i < j \leq 4 }(|A_i\cap A_j|+|B_i \cap
B_j|).
\end{eqnarray*}
Combining the equation above with equation (\ref{eq:ab0}),  we have
\begin{equation}\label{eq:4}
4u_0+ u_1+u_3+ 4u_4=-2n^2+\sum_{1\leq i < j \leq 4 }(|A_i\cap
A_j|+|B_i \cap B_j|).
\end{equation}
Lemma \ref{l1} implies
 $|A_i \cap A_j| \geq (\alpha n)^2$ and $|B_i \cap B_j|
\geq (n-\alpha n)^2$ for $(i,j) \in \{(1,2),(2,3),(3,4),(1,4)\}$. We
get
\begin{eqnarray}
u_1+u_3+4u_0+4u_4& \geq & 2n^2-8\alpha n^2+8\alpha^2 n^2+|A_1\cap
A_3|\nonumber\\
 &+& |A_2\cap A_4| +|B_1 \cap B_3|+|B_2\cap B_4|. \label{eq:6}
\end{eqnarray}
Let $E$ be the collection of all even-colored 4-term progressions
and $O$ be the collection of all odd-colored 4-term progressions. We
have $|E|=u_0+u_2+u_4$ and $|O|=u_1+u_3$.  Inequality (\ref{eq:6})
together with $\sum_{i=0}^4 u_i =n^2$ give that
\begin{eqnarray}
 m_4(\Z_n,c)&=&u_0+u_4 \nonumber \\
&=& \frac{1}{4}(u_1+u_3+4u_0+4u_4+|E|-n^2) \nonumber\\
& \geq& (\frac{1}{4}-2\alpha+2\alpha^2)n^2+ \frac{|E|}{4}+
\frac{1}{4}(|A_1\cap A_3| +|B_1 \cap B_3|)
\nonumber\\
&+& \frac{1}{4}(|A_2 \cap A_4|+|B_2 \cap
B_4|).\label{eq:7}
\end{eqnarray}
We aim to modify the method in \cite{wolf} to find a lower bound on
$|E|$ which gives a lower bound on $m_4(\Z_n,c)$. Assume $S$ is a
3-term progression in $\Z_n$.  Let $p_S$ be the number of
even-colored 4-APs containing $S$ and $q_S$ be the number of
odd-colored 4-APs containing $S$.  Observe that $p_S+q_S=2$. If
$a,a+d,a+2d$ is a 3-AP, then it determines a pair of integers $x,y
\in \Z_n$ such that $x,a,a+d,a+2d$ and $a,a+d,a+2d,y$ are two 4-APs
containing $S$; the pair $(x,y)$ is the {\it frame pair} of $S$.  We
have
$$
\mathbb E_S p_S=2|E| \ \mbox{and} \  \mathbb E_S q_S=2|O|,
$$
where the expectation operator $\mathbb E_S$ runs over all 3-APs.
The following equality which ensures us to obtain a  lower bound on
$E$. We have
\begin{eqnarray}
2|E|&=&2|O|+ \E_S(p_S-q_S) \nonumber\\
&=& 2(n^2-|E|)-\E_S(|p_S-q_S|)+2\E_S(p_S-q_S|p_S > q_S).
\label{eq:9}
\end{eqnarray}
Solving for $|E|$ in equation (\ref{eq:9}) gives
\begin{equation}
|E|=\frac{1}{2} n^2-\frac{1}{4}\E_S(|p_S-q_S|)+
\frac{1}{2}\E_S(p_S-q_S|p_S > q_S).
\end{equation}
We have the following claim which will be proved at the end of this
section.
\begin{claim} \label{claim2}
$\E_S(p_S-q_S|p_S > q_S) \geq n^2/12$ for any positive integer
$n$.
\end{claim}
Observe that $|p_S-q_S| \not = 0$ if and only if the frame pair of
$S$ is monochromatic. Furthermore, $|p_S-q_S|=2$ if $|p_S-q_S| \not
= 0$. Note that when $n$ is prime,  each frame pair belongs to a
unique 3-term progression as 4 is invertible in $\Z_n$. However, if
$n$ is not prime, then each frame pair may belong to more than one
3-term progression or does not belong to any 3-term progressions. We
will compute the value of $E_S(|p_S -q_S|)$ case by case according
to $n$ modulo $4$.

{\bf Case 1:} $n \equiv 1, 3 \mod 4$. In this
case, each frame pair belongs to a unique 3-term progression since 4
is invertible in $\Z_n$. We have $\E_S(|p_S-q_s|)$ equals twice of
the number of monochromatic pairs in the coloring $c$, that is
$\E_S(|p_S-q_s|)=2(\alpha n)^2+2(n-\alpha n)^2$. We obtain
 $$
  |E| \geq \alpha(1-\alpha)n^2+\frac{n^2}{24}.
$$
By Lemma \ref{l1}, we have $|A_1\cap A_3|=|A_2 \cap A_4| \geq
(\alpha n)^2$ and  $|B_1\cap B_3|=|B_2 \cap B_4| \geq (n-\alpha
n)^2$. Therefore, in this case, inequality (\ref{eq:7}) is
\begin{equation} \label{eq:11}
  m_4(\Z_n,c) \geq \frac{(3-11\alpha-11
\alpha^2)n^2}{4}+\frac{n^2}{96}.
\end{equation}
It is straightforward to check that the minimum value of  the right
hand side of inequality (\ref{eq:11}) is $7n^2/96$ and it is
achieved at $\alpha=1/2$. We have $m_4(\Z_n) \geq  7/96$ in
this case.

{\bf Case 2:} $n \equiv 2 \mod 4$. For $0 \leq i \leq
3$, let $\Z_n^i=\{z\in \Z_n\colon z\equiv i\mod 4\}$,
$a_i=|A\cap \Z_n^i|$, and $b_i=|B\cap \Z_n^i|$.
A pair $(x,y)$ is a frame pair if and only if $y-x =4d$ for some
$d\in \Z_n$. Assume $n=4r+2$. If $d$ is a solution for $4d=y-x$,
then $d+ 2r+1$ is another solution. We have
\begin{equation}\label{eq:12}
\E_S(|p_S-q_S|)=4(a_0+a_2)^2+4(a_1+a_3)^2+4(b_0+b_2)^2+4(b_1+b_3)^2.
\end{equation}
By the same argument, we have
$$
|A_1 \cap A_3|=|A_2 \cap A_4|=2(a_0+a_2)^2+2(a_1+a_3)^2
$$
and
$$
|B_1 \cap B_3|=|B_2 \cap B_4|=2(b_0+b_2)^2+2(b_1+b_3)^2.
$$
Therefore, $|E| + |A_1\cap A_3|+|A_2 \cap A_4|+|B_1 \cap B_3|+|B_2
\cap B_4|$ is at least
$$\frac{1}{2}n^2+3((a_0+a_2)^2+(a_1+a_3)^2+(b_0+b_2)^2+(b_1+b_3)^2)
+\frac{n^2}{24}.$$
We have the following inequality
\begin{eqnarray}
|E| &+& |A_1\cap A_3|+|A_2 \cap A_4|+|B_1 \cap B_3|+|B_2 \cap B_4| \label{eq:13}\\
  & \geq & \frac{1}{2}n^2+\frac{3}{2}(\sum_{i=0}^3
  a_i)^2+\frac{3}{2} (\sum_{i=0}^3b_i)^2 +\frac{n^2}{24}.\nonumber\\
  &=&\left( \frac{1}{2} + \frac{3}{2}\alpha^2+  \frac{3}{2}(1-\alpha)^2+
\frac{1}{24} \right)n^2. \nonumber
\end{eqnarray}
Combining inequalities (\ref{eq:7}) and  (\ref{eq:13}), we get
$$
m_4(\Z_n,c) \geq
\frac{(3-11\alpha+11\alpha^2)n^2}{4}+\frac{n^2}{96}.
$$
Note the minimum is reached at $\alpha=1/2$. It follows
$m_4(\Z_n)\geq 7/96$.

{\bf Case 3:} $n \equiv 0 \mod 4$. This method
 fails in this case; which suggests that it is possible to
find a good 2-coloring of $\mathbb Z_n$ which contains  few
monochromatic 4-term progressions  when $n$ is a multiple of $4$.
Replacing the terms on the right hand side of inequality
(\ref{eq:6}) by the lower bounds from Lemma \ref{l1}, we
obtain
\begin{equation} \label{eq:14}
u_1+u_3+4u_0+4u_4 \geq 4n^2-12\alpha n^2+12 \alpha^2 n^2.
\end{equation}
Combining with $\sum_{i=0}^4 u_i=n^2$, we have
\begin{equation}\label{eq:15}
u_0+u_4 \geq \frac{u_2}{3}+1-4\alpha n^2+4\alpha^2 n^2 \geq
\frac{u_2}{3}.
\end{equation}
The remark following the proof of Theorem 4.4 in  \cite{cameron}
gives
\begin{equation}\label{eq:16}
u_0+u_2+u_4 \geq \frac{8n^2}{33}.
\end{equation}
Combining inequalities (\ref{eq:15}) and (\ref{eq:16}), we get
$$
  m_4(\Z_n,c)=u_0+u_4 \geq \frac{2n^2}{33}.
$$
It implies $m_4(\Z_n) \geq 2/33$. We completed the proof of
Theorem \ref{t:m4lb}. \hfill $\square$

We finish this section by proving Claim 1.

 {\bf Proof of Claim 1:} Observe that $p_S > q_S$ if and only if the
coloring pattern of the $5$-APs ($S$ and its frame
pair $(x,y)$) is in  the following set

$
F=\{(1,1,1,1,1),(1,0,0,1,1),(1,0,1,0,1),(1,1,0,0,1),\\
\hspace*{1.5cm}(0,0,0,0,0),(0,1,1,0,0),(0,1,0,1,0),(0,0,1,1,0) \}. $

Moreover, for each $S$, $p_S-q_S=2$ if $p_S > q_S$. Therefore the
value of  $\E_S(p_S-q_S|p_S > q_S)$ is twice of the number of
increasing 5-term progressions with coloring pattern from $F$. Using
an exhaustive search, one can show that for any 2-coloring of
$[46]$, there is at least one increasing 5-AP of coloring pattern in
$F$.


A further computation shows that any $2$-coloring of $[74]$ contains
at least 27 increasing $5$-APs of coloring
pattern in $F$.
Note that the number of increasing
5-APs in $[74]$  with $d=1$ is $70$, the number of 5-APs  in  $[74]$
with $d=2$ is $66$, etc. The  number of $5$-APs in $[74]$ is
$$70+66+62+\cdots+6 +2=648.$$
For any $2$-coloring of $\Z_n$, the number of $74$-APs
is exactly $n^2$; each of them (degenerated or not)
 contains $27$  $5$-APs of coloring
pattern in $F$. Each $5$-AP with coloring pattern in $F$  is counted
at most $648$-times. Thus, the number of  5-APs with coloring
pattern in $F$ is at least
$$\frac{27}{648}n^2=\frac{1}{24}n^2.$$
 Thus  we have
$$\E_S(p_S-q_S|p_S > q_S) \geq \frac{n^2}{12}.$$
We finished the proof of the claim. \hfill $\square$


\begin{thebibliography}{9}
\bibitem{BCG}S.~Butler, K.~Costello and R.~Graham. Finding patterns avoiding many
monochromatic constellations.  {\it Experiment. Math}. {\bf 19}
(2010), 399-411.
\bibitem{cameron}P.~Cameron, J.~Cilleruelo and O.~Serra. On monochromatic solutions of
equations in groups. {\it Rev. Mat. Iberoam}.  {\bf 23}(2007),
385-395.

\bibitem{da} B.~Datskovsky. On the number of monochromatic Schur triples. {\it Adv.
Appl. Math}. {\bf 31}  (2003),193-198.
\bibitem{fgr} P.~Frankl, R.~Graham and V.~R\"odl.
Quantitative theorems for regular systems of equations. {\it J.
Combin. Theory Ser. A}. {\bf 47} (1988), 246-261.

\bibitem{prs} P.~Parillo, A.~Robertson and D.~Saracino. On the asymptotic
minimum number of monochromatic 3-term arithmetic progressions. {\it
J. Combin.Theory Ser. A}. {\bf 115} (2008), 185-192.

\bibitem{pick} G.~Pick. Geometrisches zur Zahlenlehre
{\it Naturwiss. Z. Lotus (Prag)}. (1899), pp. 311-319.

\bibitem{RZ} A.~Robertson and D.~Zeilberger.
A 2-coloring of $[1,N]$ can have $(1/22)N^2+O(N)$ monochromatic
Schur triples, but not less!. {\it Elec. J. Combin}. {\bf 5} (1998),
Research Paper 19.
\bibitem{schoen} T.~Schoen. The number of monochromatic Schur
triples. {\it Europ. J. Combinatorics}. {\bf 20} (1999), 855- 866.

\bibitem{thanatipanonda} T.~Thanatipanonda.
On the monochromatic Schur triples type problem. {\it Elec. J.
Combin.} {\bf 16} (2009), Research Paper 14.

\bibitem{thomason} A.~Thomason.
A disproof of a conjecture of Erd\H{o}s in Ramsey theory. {\it J.
London Math. Soc.} {\bf 39} (1989), 246-255.

\bibitem{waerden}  Van~der~Waerden. Beweis einer Baudetschen
Vermutung.
{\it Nieuw Arch. Wisk.} {\bf 15} (1927),
212-216.

\bibitem{wolf} J.~Wolf. The minimum number of monochromatic 4-term
 progressions in $\mathbb Z_p$. {\it J.of.Combinatorics}. {\bf 1}
(2010), 53-68.

\end{thebibliography}
\end{document}